\documentclass[12pt,a4paper]{article}
\usepackage[cp1251]{inputenc}
\usepackage{graphicx}
\usepackage{amssymb}
\usepackage{amsmath}
\usepackage{amsthm}
\usepackage{euscript}
\numberwithin{equation}{section}
\newtheorem{definition}{Definition}
\newtheorem{theorem}{Theorem}
\newtheorem{lemma}{Lemma}
\newtheorem{conjecture}{Conjecture}
\newtheorem{rem}{Remark}

\newenvironment{remark}{\begin{rem}\rm}{\end{rem}}

\newtheorem{ex}{Example}

\newcommand{\RP}[0]{\mathbb{RP}}

\newcommand{\R}[0]{\mathbb{R}}

\newcommand{\Z}[0]{\mathbb{Z}}

\newcommand{\im}[0]{\mathrm{Im}}

\newcommand{\be}[0]{\begin{equation}}
\newcommand{\ee}[0]{\end{equation}}
\newcommand{\bez}[0]{\begin{equation*}}
\newcommand{\eez}[0]{\end{equation*}}
\newcommand{\bl}[0]{\begin{lemma}}
\newcommand{\el}[0]{\end{lemma}}

\newcommand{\ep}[0]{$\hspace{\fill} \square$}
\newcommand{\paragraf}[1]{\par
\bigskip{\centerline{\bf #1}}\medskip}

\newcommand{\abs}[1]{\begin{quotation} {\small
 \centerline{{\bf Abstract}}\smallskip
#1}
\end{quotation}}
\newcommand{\prove}[0]{
\smallskip \noindent{\it Proof.}~}

\author{I.\,N.~Shnurnikov\footnote{Yaroslavl State University, Delone Laboratory of Discrete and Computational Geometry. This work is supported by the Russian Government project 11.G34.31.0053.}}
\title{On the number of connected components in complements to arrangements of submanifolds}
\date{}
\begin{document}
\maketitle
\abs{
 We consider arrangements of $n$ connected codimensional one submanifolds in closed $d$ -- dimensional manifold $M$. Let $f$ be the number of connected components of the complement in $M$ to the union of submanifolds. We prove the sharp lower bound for $f$ via $n$ and homology group $H_{d-1}(M)$. The sets of all possible $f$~--~values for given $n$ are studied for hyperplane arrangements in real projective spaces and for subtori arrangements in $d$ -- dimensional tori.
}

\paragraf{Introduction}
The theory of plane arrangements in affine or projective spaces has been investigated rather thoroughly, see the book of P.~Orlic, H.~Terao, \cite{Orlic_1992} and V.\,A.~Vassiliev's review \cite{Vasil'ev}. Inspired by a conjecture of B.~Gr\"unbaum \cite{Grunbaum 72}, N.~Martinov \cite{Martinov 93} found all possible pairs $(n,f)$ such that there is a real projective plane arrangement of $n$ pseudolines and $f$ regions. It turns out, that some facts concerning arrangements of hyperplanes or oriented matroids could be generalized to arrangements of submanifolds, see P.~Deshpande dissertation \cite{Deshpande}. So we are going to study the sets $F(M,n)$ of connected components numbers of the complements in the closed manifold $M$ to the unions of $n$ closed connected codimensional one submanifolds. Sometimes it seems reasonable to restrict the type of submanifolds, for example, author \cite{Shnurnikov 11} found sets $F(M,n)$ of region numbers in arrangements of $n$ closed geodesics in the two dimensional torus and the Klein bottle with locally flat metrics.

\paragraf{Homological bound of the number of connected components}

Let $M^n$ be connected $n$--dimensional smooth compact manifold without boun\-dary, let $A_i\subset M^n$ be distinct connected $(n-1)$--dimensional closed submanifolds in $M^n$ for $1 \leq i \leq k$. Let us consider the union
$$
A=\bigcup_{i=1}^kA_i.
$$
We shall denote by $f$ the number $|\pi_0(M^n\setminus A)|$ of connected components of the complement to $A$ in $M^n$. Let $UA$ be regular open neighbourhood of $A$ in $M^n$. Let
\bez
M^n \setminus UA \cong \bigsqcup_{j=1}^f N_j,
\eez
where $N_j$ are the connected components of the complement to $UA$ in $M^n$. If $M^n$ and all submanifolds $A_i$ are orientable, then we assume $G=\Z$. If some $A_i$ or $M^n$ is not orientable, then $G=\Z_2$.

\bl
\label{lemma homology H >= k}
If closed $(n-1)$ -- dimensional submanifolds $A_i \subset M^n$, $i=1, \dots, k$ intersect each other transversally, then
\bez
\dim H_{n-1}(UA,G)=\dim H_{n-1}(A,G) \geq k
\eez
\el

\prove
The regular neighbourhood of $UA$ is homotopically equivalent to $A$ and so all homology groups of $A$ and $UA$ are the same.
By induction on $k$ let us prove
\bez
\dim H_{n-1}\left(\cup_{i=1}^k A_i, G\right)\geq k.
\eez
It is obvious for $k=1$ because for connected closed $(n-1)$ -- dimensional manifold $H_{n-1}(A_1,G)\cong G$. Suppose the statement is true for $k-1$ submanifolds and let us prove it for $k$ submanifolds. Let
$$
A'=\bigcup_{i=1}^{k-1}A_i.
$$
Then by induction assumption
$$
\dim H_{n-1}(A', G) \geq k-1.
$$
By Meyer-Vietoris exact sequence for pair $A', A_k$ we have:
\bez
\longrightarrow H_{n-1}(A'\cap A_k)\longrightarrow H_{n-1}(A')\oplus H_{n-1}(A_k) \longrightarrow H_{n-1}(A'\cup A_k) \longrightarrow
\eez
As submanifolds $A_i$ and $A_j$ intersect transversally then $A'\cap A_k$ is a finite union of at most $(n-2)$ -- dimensional submanifolds in $M^n$. Hence $H_{n-1}(A'\cap A_k)=0$ and the map
\bez
H_{n-1}(A')\oplus H_{n-1}(A_k) \longrightarrow H_{n-1}(A'\cup A_k)
\eez
is monomorphic. Therefore,
\bez
\dim H_{n-1}(A'\cup A_k)\geq \dim H_{n-1}(A')+\dim H_{n-1}(A_k) \ \geq \ k.
\eez
\ep

\bl
\label{lemma homology G^f}
\bez
H_n(M^n, UA, G)\cong G^f
\eez
\el
\prove
\begin{multline}
H_n(M^n, UA, G)=\widetilde{H}_n(M^n/UA, G)= \notag\\
=\widetilde{H}_n\left(\sqcup_{j=1}^f N_j/ \sqcup_{j=1}^f \partial N_j, \ G\right)= \widetilde{H}_n\left(\vee_{j=1}^f N_j/\partial N_j, G\right)= \\
=\bigoplus_{j=1}^f  \widetilde{H}_n(N_j/\partial N_j, G)=G^f,
\end{multline}
where $n \geq 1,$ $\vee$ is one point union, $\widetilde{H}_n$ is the reduced homology group, $\partial N_j$ is the boundary of $N_j$. \ep

\begin{theorem}
\label{theorem homology f >=n-dim H}
Let $A_1, \dots, A_k$ be connected closed codimensional one submanifolds in a connected closed  $n$~--~dimensional manifold $M^n$. Suppose that the submanifolds $A_i$ intersect each other transversally and $A=\cup_i A_i$. Then
\bez
|\pi_0(M^n\setminus A)|\geq k+1 - \dim H_{n-1}(M^n,G),
\eez
where $G$ is chosen as before.
\end{theorem}
\prove
Let us write the exact homological pair sequence for inclusion $i: UA \to M^n$ with coefficients in $G$:
\bez
H_n(UA) \longrightarrow H_n(M^n) \longrightarrow H_n(M^n, UA) \longrightarrow H_{n-1}(UA) \longrightarrow H_{n-1}(M^n)\longrightarrow
\eez
Notice that
\bez
H_n(UA)=H_n(A) =0, \quad H_n(M_n)=G.
\eez
It follows from the exactness of sequence in $H_n(M^n)$, that the map
\bez
H_n(M^n) \longrightarrow H_n(M^n, UA)
\eez
is monomorphic.
By lemma \ref{lemma homology G^f}
\bez
H_n(M^n, UA, G)\cong G^f
\eez
One can see that
\bez
\dim H_{n-1}(UA) \leq \dim \im \partial_{\ast}+\dim \im i_{\ast},
\eez
where the homomorphisms are
\bez
\partial_{\ast}:H_n(M^n, UA) \longrightarrow H_{n-1}(UA), \quad \quad i_{\ast}: H_{n-1}(UA) \longrightarrow H_{n-1}(M^n).
\eez

Notice that
\bez
\dim \im i_{\ast} \leq \dim H_{n-1}(M^n), \quad \dim \im \partial_{\ast} \leq f-1.
\eez
By lemma \ref{lemma homology H >= k} $\dim H_{n-1}(UA) \geq k$ and so
$k \leq f-1 + \dim H_{n-1}(M^n).$
\ep

\begin{remark}
One can see that the inequality of the theorem is sharp for arrangements of
$$
n\geq \dim H_{n-1}(M^n)
$$
submanifolds in projective spaces, spheres, $n$~--~dimensional tori and Riemann surfaces of genus $g$.
\end{remark}

\paragraf{Toric arrangements}
\begin{definition}
By a flat $d$ -- dimensional torus $T^d$ we mean a quotient of affine $d$ -- dimensional space by a nondegenerate $d$ -- lattice $Z^d$ (which is not surely integer lattice). A codimensional one subtorus is given by equation
\bez
\sum_i a_ix_i=c,
\eez
where $a_i$ are rational, $x_i$ are coordinates of $\R^d$ in some lattice basis, $c$ is real.
\end{definition}

A codimensional one subtorus is closed submanifold of $T^d$ homeomorphic to $(d-1)$ -- dimensional torus. Let $A$ be the union of $n$ codimensional one subtori in the flat $d$~--~dimensional torus $T^d$. Consider the connected components of the complement $T^d\setminus A$; denote the number of connected components by $f=|\pi_0\left(T^d\setminus A\right)|$; let $F\left(T^d,n\right)$ be the set of all possible numbers $f$.

\begin{theorem}
\label{theorem F set for d-dimens. torus}
For $n>d$
\bez
F\left(T^d,n\right) \supseteq \{n-d+1, \dots, n\} \cup \{l \in \mathbb{N} \ | \ l \geq 2(n-d).\}
\eez
For $2\leq n \leq d$ we have $F\left(T^d,n\right)=\mathbb{N}$.
\end{theorem}
\prove
Let $T^d=\R^d/Z^d$ and $e$ be the basis of $Z^d$. Let $(x_1, \dots, x_d)$ be the coordinates of $\R^d$ in the basis $e$. We shall construct examples for $\leq n$ and $\geq 2n-2d$ regions separately.

Let us consider $n$ hyperplanes  in $\R^d$ (an equation corresponds to a hyperplane):
\begin{gather*}
x_i=0, \quad 1 \leq i \leq k,\\
x_{k+1}=c_{i-k}, \quad k+1 \leq i \leq n
\end{gather*}
for some integer $k$, $0 \leq k \leq d-1$ and real $c_{i-k}$ with different frctional parts. By the factorization map $\R^d \to \R^d/Z^d$ we shall get a set  $\left\{T_i^{d-1}, i=1,\dots,n\right\}$ of $n$ codimensional one subtori. And the complement is homeomorphic to the prime product
$$
T^d\setminus\bigcup_i T_i^{d-1}\approx \R^k\times\left(S^1\setminus\{p_1,\dots,p_{n-k}\}\right)\times\left(S^1\right)^{d-k-1},
$$
where $S^1\setminus\{p_1,\dots,p_{n-k}\}$ denotes a circle without $n-k$ points. Hence the number of complement regions equals $n-k$, for an integer $k$ such that $0 \leq k \leq d-1$.

Now let us take integer nonnegative $k$ and construct an arrangement with $2n-2d+k$ connected components of the complement. We shall determine the subtori by equations:
\begin{gather*}
x_i=0, \quad \text{for} \quad 2\leq i \leq d,\\
x_2=kx_1+ \frac12,\\
x_1=c_j \quad \text{for} \quad j=1, \dots, n-d,
\end{gather*}
whereas numbers $kc_j+\frac 12$ are not integer for any $j$. ( This means that the intersection of three subtori
$$
x_2=kx_1+ \frac12, \ x_1=c_j, \ x_2=0
$$
is an empty set.) One may see that
$$
T^d\setminus\bigcup_{i=3}^d\{x_i=0\}\approx T^2\times \R^{d-2}.
$$
In the two\,--\,dimensional torus the equations
\begin{gather*}
x_2=0,\\
x_2=kx_1+\frac 12,\\
x_1=c_j \ \text{for} \ j=1, \dots, n-d
\end{gather*}
produce the arrangement of $n-d+2$ closed geodesics. The geodesics' union divides the torus into $2n-2d+k$ connected components (for more details on arrangements of closed geodesics in the flat torus see author's paper \cite{Shnurnikov 11}).
\ep

\begin{conjecture}
It seems believable that the inclusion in the theorem is indeed the equality for all $d\geq 2$ and $n\geq d$. Yet the equality is proved for $d=2$ in \cite{Shnurnikov 11}.
\end{conjecture}

\paragraf{Sets of region's numbers in hyperplane arrangements}

By an arrangement of $n$ hyperplanes in the real projective space $\RP^d$ we mean a set of $n$ hyperplanes, such that there are no point belonging to all the hyperplanes. The arrangement produce the cell decomposition of the $\RP^d$; let $f$ denotes the number of open $d$\,--\,cells.
Let $F_n^{(d)}$ denotes the set of all possible numbers $f$ arising in arrangements of $n$ hyperplanes in $\RP^d$. Let $m$ be the maximal number of hyperplanes, passing through one point.

\begin{lemma}
\label{theorem f via n m}
For arrangements of $n$ hyperplanes in $\RP^d$ we have
\bez
f \geq (m-d+1)\sum_{j=0}^{[\frac d2]}\frac{C_n^{d-2j}}{C_{m-2j}^{d-2j}}.
\eez
\end{lemma}
\prove
It follows from Zaslavsky formula for number of regions and some inequalities concerning the M\"obius function of the arrangement poset.
  %\ref{theorem Zaslavsky} and \ref{theorem sum mu via n, m, i}.
\ep

\begin{lemma}
\label{lemma f via n m d weak}
For arrangement of $n$ hyperplanes in the real projective space $\RP^d$
$$
f\geq(n-m+1)(m-d+2)2^{d-2}.
$$
\end{lemma}
\prove
Let $m$ hyperplanes $A_1, \dots, A_m$ have nonempty intersection $Q$ ($Q$ is a point). The family $A_1, \dots, A_m$ is a cone over some arrangement $B$ of $m$ planes in $\RP^{d-1}$. The number $f(B)$ of regions in arrangement $B$ could be estimated (see Shannon paper \cite{Shannon 76}, where this result is referred to McMullen) as:
\bez
f(B)\geq (m-d+2)2^{d-2}.
\eez
Each of the remaining hyperplane of the former arrangement intersects the family $A_1, \dots, A_k $ by an arrangement $B_i$, projective equivalent to $B$. Thus
$$
f \geq f(B)+\sum_if(B_i)=(n-m+1)f(B).
$$
\ep

\begin{theorem}
\label{theorem first 4 numders}
Let $d \geq 3 $ and $n \geq 2d+5$. Then the first four increasing numbers of $F_n^{(d)}$ are the following:
\bez
(n-d+1)2^{d-1}, \quad 3(n-d)2^{d-2}, \quad (3n-3d+1)2^{d-2}, \quad 7(n-d)2^{d-3}.
\eez
\end{theorem}
\prove
We are going to prove that the four mentioned numbers are the only realizable ones among numbers not greater then $7(n-d)2^{d-3}$. After it one may see how to construct examples of arrangements with required numbers $f$. Let us prove that if $m \leq d+1$, then
 $$
f \geq  7(n-d)2^{d-3}.
$$
For $m=d$ we have an arrangement of hyperplanes in general position and the number of regions is the largest possible. If $m=d+1$, then by lemma \ref{theorem f via n m} we have
\bez
f\geq \frac {C_{n}^{d+1}}{n-d}= \frac{n}{3}\frac{(n-1)}{d+1}\frac{(n-2)}{d}\frac{(n-3)}{d-1}\dots \frac{(n-d+2)}{4}(n-d+1)\geq 7\cdot2^{d-3}(n-d)
\eez
because $n\geq 2d+5.$

Now we prove the theorem for $d=3,$ $n\geq 11$.
Let us consider three cases.

{\it 1.} If $m=n-1$, then $f=2\varphi$, where $\varphi \in F_{n-1}^{(2)}$. The set $F^{(2)}_{n-1}$ is known due to N.~Martinov \cite{Martinov 93}
\bez
\{f \in F^{(2)}_{n-1} \ |\  f \leq 4n-16 \}= \{2n-4, 3n-9, 3n-8, 4n-16\}.
\eez

{\it 2.} $m=n-2$. The arguments are the same as in the inductive step further (Martinov theorem \cite{Martinov 93} for the set $F_n^{(2)}$ is also used).

{\it 3.} If $5 \leq m \leq n-3$, then by using lemma \ref{lemma f via n m d weak} we have
\bez
f\geq 2(n-m+1)(m-1) \geq 8n-32 \geq 7n-21
\eez
for $n \geq 11.$

Now we use induction on $d\geq 3$. Base is the validity of the theorem for $d=3$. The assumption is the validity of the theorem for all integers $3\leq d' < d$ and $n' \geq 2d'+5$. To prove the induction step we shall consider three cases.

{\it 1.} If $m=n-1$, then $f=2\varphi$, where $\varphi \in F_{n-1}^{(d-1)}$. By induction assumption for the set $F_{n-1}^{(d-1)}$ (note that $n-1\geq 2(d-1)+1$) we get that either $\varphi$ is equal to one of four numbers
\bez
(n-d+1)2^{d-2}, \quad 3(n-d)2^{d-3}, \quad (3n-3d+1)2^{d-3}, \quad 7(n-d)2^{d-4},
\eez
or $\varphi >7(n-d)2^{d-4}.$

{\it 2.} $m=n-2$. Consider $n-2$ hyperplanes $p_1, \dots, p_{n-2}$, passing through one point. These hyperplanes cut $\RP^d$ into $\varphi$ regions and $\varphi \in F^{(d-1)}_{n-2}$. Let $l$ denote the intersection of the two remaining hyperplanes. By the inductive assumption we have either
$$
\varphi=(n-d)2^{d-2} \quad \text{or} \quad \varphi \geq 3(n-d-1)2^{d-3}
$$
(note that assumption may be used as $n-2\geq 2(d-1)+5$).
If
$$
l\in \bigcup_{i=1}^{n-2}p_i
$$
then $f=3\varphi$ and the case is over.
If
$$
l\notin \bigcup_{i=1}^{n-2}p_i
$$
then let $B$ be the set of planes $p_i\cap l$ in the $l$, where $l$ is regarded as the ambient $(d-2)$\,--\,dimensional projective space. One may prove, that $B$ is an arrangement of at least $n-3$ planes in $l$. Then $f(B)\geq (n-d)2^{d-3}$ by Shannon theorem \cite{Shannon 76}. Since
$$
f=3\varphi+f(B)\geq 7(n-d)2{d-3},
$$
the case is over.

{\it 3.} If $d+2 \leq m \leq n-3$ then by lemma \ref{lemma f via n m d weak} we have
\bez
f\geq(n-m+1)(m-d+2)2^{d-2} \geq (4n-4d-4)2^{d-2} \geq 7(n-d)2^{d-3}
\eez
for $n \geq d+8.$
\ep

\begin{lemma}
\label{lemma f via n m d quadratic}
For arrangement of $n$ hyperplanes in the real projective space $\RP^d$
\bez
f \geq 2\frac{n^2-n}{m-d+5}.
\eez
\end{lemma}
\proof
It follows from the similar inequality for arrangement of lines in the projective plane, see details in \cite{Shnurnikov 10}.
\ep

\begin{theorem}
\label{theorem first 36 numbers}
First 36 increasing numbers of the set $F^{(3)}_n$ for $n\geq 50$ are the following (i.e. all realizable numbers up to $12n-60$)
\begin{gather*}
4n-8, \ 6n-18, \ 6n-16, \ 7n-21,
\ 7n-20, \ 8n-32, \ 8n-30, \ 8n-28,\\
8n-26, \ 9n-36, \ 9n-33, \ 9n-31,
\ 9n-30, \ 10n-50, \ 10n- 48, \ 10n-46, \\
10n- 44, \ 10n-42, \ 10n-40, \ 10n-39,\ 10n-38, \ 10n-37, \ 10n-36, \ 10n-35, \\
11n-44, \ 11n-43, \ 11n-42, \ 11n-41, \ 11n-40, \ 12n-72, \ 12n-70, \ 12n-68, \\
\ 12n-66, \ 12n-64, \ 12n-62, \ 12n-60.
\end{gather*}
\end{theorem}

\prove Let $m$ be the maximal number of hyperplanes, passing through one point. Examples for this numbers could be constructed for arrangements with $m\geq n-5$. Let us prove that there are no other realizable numbers, smaller then $12n-60$. Consider three cases.

{\it 1.} It $m\geq n-5$, then by enumeration of possibilities we have that either $f$ belongs to given set or $f \geq 12n-60.$

{\it 2.} If $8 \leq m \leq n-6$, then by lemma \ref{lemma f via n m d weak} we have
$f\geq 7n-49.$

{\it 3.} If $m \leq 7$ then by lemma \ref{lemma f via n m d quadratic}
\bez
f \geq 2\frac{n^2-n}9 \geq 12n-60
\eez
for $n\geq 50$.
\ep

\end{document}